\documentclass[12pt,twoside]{amsart}
\usepackage{amssymb}
\usepackage{amscd}
\usepackage[abbrev,alphabetic]{amsrefs}
\usepackage{hyperref}
\usepackage{here}
\usepackage{comment}
\usepackage[margin=1.25in]{geometry}
\usepackage{color}

\usepackage{fancybox}
\usepackage{ascmac}

\title[Boundedness of regular del Pezzo surfaces]
{Boundedness of regular del Pezzo surfaces over imperfect fields} 
\author{Hiromu Tanaka} 
\subjclass[2010]{14G17, 14D06.}
\keywords{imperfect fields, del Pezzo surfaces, positive characteristic.}
\address{Graduate School of Mathematical Sciences, 
The University of Tokyo, 
3-8-1 Komaba, Meguro-ku, Tokyo 153-8914, JAPAN} 
\email{tanaka@ms.u-tokyo.ac.jp}
\newcommand{\Chow}[0]{{\operatorname{Chow}}}
\newcommand{\Hilb}[0]{{\operatorname{Hilb}}}
\newcommand{\Univ}[0]{{\operatorname{Univ}}}

\newcommand{\red}[0]{{\operatorname{red}}}

\newcommand{\Proj}[0]{{\operatorname{Proj}}}
\newcommand{\Spec}[0]{{\operatorname{Spec}}}

\newtheorem{thm}{Theorem}[section]
\newtheorem{lem}[thm]{Lemma}
\newtheorem{cor}[thm]{Corollary}
\newtheorem{prop}[thm]{Proposition}
\newtheorem{claim}[thm]{Claim}

\theoremstyle{definition}
\newtheorem{ex}[thm]{Example}

\newtheorem{dfn}[thm]{Definition}

\newtheorem{rem}[thm]{Remark}

\makeatletter
  
  \@addtoreset{equation}{thm}
  \makeatother

\newcommand{\MO}{\mathcal{O}}

\newcommand{\Q}{\mathbb{Q}}
\newcommand{\Z}{\mathbb{Z}}
\newcommand{\F}{\mathbb{F}}

\newcommand{\m}{\mathfrak{m}}

\begin{document}

\maketitle

\begin{abstract}
For a regular del Pezzo surface $X$, we prove that $|-12K_X|$ is very ample. Furthermore, we also give an explicit upper bound for the volume $K_X^2$ which depends only on $[k : k^p]$ for the base field $k$. As a consequence, we obtain the boundedness of geometrically integral regular del Pezzo surfaces. 
\end{abstract}

\tableofcontents

\section{Introduction}

One of important classes of algebraic varieties {is the one of Fano varieties}. 
For example, classification of Fano varieties 
has been an interesting problem in algebraic geometry. 
Indeed, Fano varieties are classified in dimension at most three (cf. \cite{MM83}). 
Although it seems to be difficult to obtain 
complete classification in higher dimension, 
it turns out that Fano varieties form bounded families 
when we fix the dimension \cite{KMM92}. 
Apart from the boundedness, 
Fano varieties satisfy various prominent properties, 
e.g. they are rationally connected (\cite{Cam92}, \cite{KMM92}) and 
have no non-trivial torsion line bundles.

The main topic of this article is to study regular del Pezzo surfaces over imperfect fields. 
We naturally encounter such surfaces when we study minimal model program over algebraically closed fields of positive characteristic. 
The minimal model conjecture predicts that 
an arbitrary algebraic variety is birational to either a minimal model 
or a Mori fibre space $\pi:V \to B$. 
Although general fibres of $\pi$ might have bad singularities in positive characteristic 
(e.g. they are non-reduced if $\pi:V \to B$ is a wild conic bundle \cite{MS03}), 
the generic fibre $X:=V \times_B \Spec\,K(B)$ 
of $\pi$ allows only terminal singularities. 
Note that the base field $K(B)$ of $X$ is no longer a perfect field 
in general. 
Furthermore, if $\dim X=2$, then $X$ is a regular del Pezzo surface over $K(B)$.

The purpose of this article is to establish 
results related to boundedness of regular del Pezzo surfaces. 
The main results are the following two theorems. 

\begin{thm}[Theorem \ref{t-dP-va2}]\label{intro-dP-va2}
Let $k$ be a field of characteristic $p>0$. 
Let $X$ be a regular projective surface over $k$ such that 
$-K_X$ is ample and $H^0(X, \MO_X)=k$. 
Then the complete linear system $|-12K_X|$ is very ample over $k$, 
i.e. it induces a closed immersion to $\mathbb P^N_k$ for $N:=\dim_k H^0(X, \MO_X(-12K_X))-1$. 
\end{thm}

\begin{thm}[Corollary \ref{c-bd-not-gn2}, Theorem \ref{t-vol-bound2}]\label{intro-vol-bound2}
Let $k$ be a field of characteristic $p>0$.  
Let $X$ be a regular projective surface over $k$ such that 
$-K_X$ is ample and $H^0(X, \MO_X)=k$. 
Then the following hold. 
\begin{enumerate}
\item 
If $r:=\log_p[k:k^p]<\infty$, then $K_X^2 \leq \max \{9, 2^{2r+1}\}$. 
\item 
If $X$ is geometrically reduced over $k$, then $K_X^2 \leq 9$. 
\end{enumerate}
\end{thm}

\begin{rem}
Let $\mathbb F$ be an algebraically closed field of characteristic $p>0$. 
Let $\pi:V \to B$ be a Mori fibre space between normal varieties over $\F$. 
Then $V$ has at worst terminal singularities. 
Assume that $X:=V \times_B \Spec\,K(B)$ is of dimension two. 
It holds that $X$ is a regular projective surface over $k:=K(B)$ such that 
$-K_X$ is ample and $H^0(X, \MO_X)=k$. 
In this case, it holds that 
\[
r = \log_p[k:k^p]=\log_p[K(B):K(B)^p]=\dim B. 
\]
Hence, $r$ is the dimension of the base of the Mori fibre space.
\end{rem}

As a consequence, we obtain the boundedness of 
geometrically integral regular del Pezzo surfaces. 

\begin{thm}[Theorem \ref{t-dP-bdd3}]\label{intro-dP-bdd3}
There exists a flat projective morphism $\rho:V \to S$ 
of quasi-projective $\Z$-schemes 
which satisfies the following property: 
if $k$ is a field and 
$X$ is a regular projective surface over $k$ 
such that $-K_X$ is ample, $H^0(X, \MO_X)=k$, and $X$ is geometrically reduced over $k$, 
then there exists a cartesian diagram of schemes: 
\[
\begin{CD}
X @>>> V\\
@VV\alpha V @VV\rho V\\
\Spec\,k @>>> S, 
\end{CD}
\]
where $\alpha$ denotes the induced morphism. 
\end{thm}

\begin{rem}
We fix a field $k$ such that $[k:k^p]<\infty$. 
Then Theorem \ref{intro-dP-va2} and Theorem \ref{intro-vol-bound2}(1)  
show that 
if $X$ is a regular projective surface over $k$ such that 
$-K_X$ is ample and $H^0(X, \MO_X)=k$, 
then $K_X^2$ is bounded and $|-12 K_X|$ is very ample. 
It is tempting to conclude the boundedness of these surfaces. 
However, we obtain the boundedness only for 
the geometrically reduced case as in Theorem \ref{intro-dP-bdd3}. 
In our proof, we use the following two facts (cf. Proposition \ref{p-chow-var}): 
\begin{enumerate}
\item 
A Chow variety  is a coarse moduli space 
(cf. \cite[Ch. I, Section 3, Section 4]{Kol96}), 
which does not have enough information on non-geometric points. 
\item 
The proof of the inequality $\deg X \geq 1+ {\rm codim}\,X$ for 
nondegenerate varieties $X \subset \mathbb P^N$ 
(cf. \cite[Proposiiton 0]{EH87}) works for varieties only over algebraically closed fields. 
\end{enumerate}
\end{rem}

Theorem \ref{intro-dP-bdd3} immediately implies the following corollary. 
{Indeed, Theorem \ref{intro-dP-bdd3} establishes the equivalence between 
the boundedness of $\dim_k H^1(X, \MO_X)$ and 
the boundedness of geometrically intergal regular del Pezzo surfaces.}

\begin{cor}[Corollary \ref{c-dP-bdd3}]\label{intro-dP-bdd4}
There exists a positive integer $h$ 
which satisfies the following property: 
if $k$ is a field of characteristic $p>0$ and 
$X$ is a regular projective surface over $k$ 
such that $-K_X$ is ample, $H^0(X, \MO_X)=k$, and $X$ is geometrically reduced over $k$, 
then $\dim_k H^1(X, \MO_X) \leq h$. 
\end{cor}

The original motivation of the author was  
to establish results toward 
the Borisov--Alexeev--Borisov (BAB, for short) conjecture for threefolds 
over algebraically closed fields of positive characteristic. 
One of the steps of the proof of BAB conjecture in characteristic zero 
is to apply induction on dimension by using Mori fibre spaces 
(\cite{Bir1}, \cite{Bir2}). 
If we adopt a similar strategy for threefolds in positive characteristic, 
it is inevitable to treat 
three-dimensional del Pezzo fibrations. 
In characteristic zero, we may apply the induction hypothesis 
for general fibres, 
whilst we probably need to treat generic fibres in positive characteristic 
as replacements of general fibres. 
Thus, the author originally wanted to prove 
the boundedness of geometrically integral $\epsilon$-klt log del Pezzo surfaces. 
Although Theorem \ref{intro-dP-bdd3} is weaker than this goal, 
the author hopes that our results and techniques will be useful 
to establish such generalisation.

\subsection{Description of proofs}

\subsubsection{Sketch of Theorem \ref{intro-dP-va2}}

Let $k$ be a field of characteristic $p>0$. 
Let $X$ be a regular projective surface over $k$ such that 
$-K_X$ is ample and $H^0(X, \MO_X)=k$. 
Let us overview how to find a constant $m>0$ such that $|-mK_X|$ is very ample. 
Combining known results, it is not difficult to show that 
$|-nK_X|$ is base point free for some constant $n>0$ (cf. the proof of Theorem \ref{t-dP-va1}). 
Then the problem is reduced to show the following theorem of Fujita type.

\begin{thm}[Theorem \ref{t-va-criterion}]\label{t-va-criterion}
Let $k$ be a field of characteristic $p>0$. 
Let $X$ be a $d$-dimensional regular projective variety over $k$. 
Let $A$ be an ample invertible sheaf on $X$ and 
let $H$ be an ample globally generated invertible sheaf on $X$. 
Then $\omega_X \otimes_{\MO_X} H^{d+1} \otimes_{\MO_X} A$ is very ample over $k$. 
\end{thm}

Indeed, by applying this theorem for $A:=\MO_X(-K_X)$, $H:=\MO_X(-nK_X)$, 
and $m:=3n$, 
it holds that $|-mK_X|$ is very ample. 
We now give a sketch of the proof of Theorem \ref{t-va-criterion}. 
Note that Theorem \ref{t-va-criterion} is known for the case when $k$ is 
an algebraically closed field (\cite[Theorem 1.1]{Kee08}). 
Thus, if $k$ is a perfect field, then we are done by taking the base change to the algebraic closure. 
However, if $k$ is an imperfect field, then 
the base change $X \times_k \overline k$ might be no longer regular. 
Hence, the problem is not directly reduced to the case when $k$ is algebraically closed. 
On the other hand, our strategy is very similar to the one of 
\cite{Kee08} and we use also the base change $X \times_k \overline k$. 

The outline is as follows. 
It is easy to reduce the problem to the case when $k$ is an $F$-finite field, i.e. 
$[k:k^p]<\infty$. 
Fix $e \in \Z_{>0}$. 
Then, for the $e$-th iterated absolute Frobenius morphism 
\[
\Phi_e:X_e \to X, \qquad X_e:=X, 
\]
the composite morphism $\beta:X_e \to X \xrightarrow{\alpha} \Spec\,k$ is of finite type, 
where $\alpha:X \to \Spec\,k$ denotes the structure morphism. 
We consider $X_e$ as a $k$-scheme via $\beta$. 
For the algebraic closure $\kappa:= \overline k$ of $k$, 
consider the base change of $\Phi_e$ by $(-) \times_k \kappa$: 
\[
\Psi_e: Y_e \to Y, \quad Y:=X \times_k \kappa, \quad Y_e:=X_e \times_k \kappa. 
\]
Since the trace map $(\Phi_e)_*\omega_{X_e} \to \omega_X$ of Frobenius is surjective, 
also the trace map $(\Psi_e)_*\omega_{Y_e} \to \omega_Y$ is surjective. 
Using Mumford's regularity, we can show that 
$(\Psi_e)_*\omega_{Y_e} \otimes H'^{d+1} \otimes A' \otimes \m_y$ is globally generated for any closed point $y$ of $Y$ and $e \gg 0$, 
where $H'$ and $A'$ are the pullbacks of $H$ and $A$, respectively. 
Then $\omega_Y \otimes H'^{d+1} \otimes A' \otimes \m_y$ is globally generated. 
Therefore,  $\omega_Y \otimes H'^{d+1} \otimes A'$ is very ample, 
hence so is $\omega_X \otimes H^{d+1} \otimes A$. 
For more details, see Section \ref{s-va}.

\subsubsection{Sketch of Theorem \ref{intro-vol-bound2}}

Both (1) and (2) of Theorem \ref{intro-vol-bound2} 
are consequences of the following theorem. 

\begin{thm}[{Corollary \ref{c-bd-not-gn2}}]\label{intro-bd-not-gn2}
Let $k$ be a field of characteristic $p>0$. 
Let $X$ be a regular del Pezzo surface over $k$ such that $H^0(X, \MO_X)=k$. 
Then the following hold {(for the definition of $\epsilon(X/k)$, see Definition \ref{d-epsilon})}. 
\begin{enumerate}
\item If $p \geq 5$, then $K_X^2 \leq 9$. 
\item If $p =3$, then $K_X^2 \leq \max \{9, 3^{\epsilon(X/k)+1}\}$. 
\item If $p =2$, then $K_X^2 \leq \max \{9, 2^{\epsilon(X/k)+3}\}$. 
\end{enumerate}
In particular, if $X$ is geometrically reduced over $k$, 
then it is known that $\epsilon(X/k)=0$, hence we obtain $K_X^2 \leq 9$. 
\end{thm}

Let us overview some of the ideas of the proof of Theorem \ref{intro-bd-not-gn2}. 
If $X$ is geometrically normal, 
then the assertion follows from a combination of known results (cf. the proof of Theorem \ref{t-bd-not-gn2}(1)). 
Hence, we only treat the case when $X$ is not geometrically normal. 
In particular, we may assume that $p \leq 3$ (cf. Theorem \ref{t-ell-bdd}(1)). 

For $Z:=(X \times_k \overline k)_{\red}^N$, 
let $g:Z \to X$ be the induced morphism. 
Then there is an effective $\Z$-divisor $D$ on $Z$ 
which satisfies the following linear equivalence (cf. Theorem \ref{t-T-PW}): 
\[
K_Z+D \sim g^*K_X.
\]
A key observation is that 
there are only finitely many possibilities for the pair $(Z, D)$ 
(Theorem \ref{t-bd-not-gn}). 
Indeed, this is enough for our purpose by the following equation 
(cf. Lemma \ref{l-epsilon-bc}): 
\[
K_X^2 = p^{\epsilon(X/k)} (K_Z+D)^2. 
\]
We now give a sketch of how to restrict the possibilities for $Z$. 
It is known that $Z$ is either a Hirzebruch surface or 
a weighted projective plane $\mathbb P(1, 1, m)$ for some $m \in \Z_{>0}$ 
(Theorem \ref{t-classify-bc}). 
For the latter case: $Z=\mathbb P(1, 1, m)$, 
it holds that $m \leq 4$ because the $\Q$-Gorenstein index is known to be bounded (Theorem \ref{t-Q-Gor}). 
Let us focus on the the case when $Z \simeq \mathbb P_{\mathbb P^1}(\MO \oplus \MO(n))$ for some $n \geq 0$. 
The goal is to prove that $n \leq 4$. 
Since $\rho(X) \leq \rho(Z)=2$, we have either $\rho(X)=1$ or $\rho(X)=2$. 
If $\rho(X)=1$, then we can show that $n=0$ 
by using Galois symmetry (Lemma \ref{l-pic-one}). 
Assume that $\rho(X)=2$. 
Then there are two extremal rays, both of which induce morphisms $X \to X'$ and $X \to X''$. 
Taking the base change to the algebraic closure, we obtain morphisms $Z \to Z'$ and $Z \to Z''$. 
The essential case is $\dim X'=2$. 
If $X'$ is not geometrically normal, then we may apply the above argument 
for $X'$, so that we deduce $n \leq 4$. 
If $X'$ is geometrically normal, then $Z'$ is canonical (Theorem \ref{t-classify-bc}), hence 
we have $n \leq 2$. 
For more details, see Section \ref{s-bdd-vol}.

\subsection{Related results}

We first review results on del Pezzo surfaces over algebraically closed fields of characteristic $p>0$. 
It is a classical result that smooth del Pezzo surfaces are classified, 
and in particular bounded. 
Then, in \cite{Ale94},  Alexeev proved the 
BAB conjecture for surfaces, 
i.e. $\epsilon$-klt log del Pezzo surfaces are bounded (cf. \cite{Jia13}). 
As for vanishing theorems, 
smooth del Pezzo surfaces over algebraically closed fields 
satisfy Kawamata--Viehweg vanishing \cite[Proposition A.1]{CT18}. 
However, if $p \in \{2, 3\}$, 
then there exist log del Pezzo surfaces violating 
Kawamata--Viehweg vanishing 
(\cite[Theorem 1.1]{Ber}, \cite[Lemma 2.4, Theorem 3.1]{CT18}, \cite[Theorem 4.2]{CT}). 
On the other hand, if $p \gg 0$, 
it is known that Kawamata--Viehweg vanishing 
holds for any log del Pezzo surfaces \cite[Theorem 1.2]{CTW17}. 
It is remarkable that this result is applied to show that 
three-dimensional klt singularities of large characteristic are 
rational singularities \cite{HW}. 

We now switch to the situation over imperfect fields. 
The first remarkable result is given by Schr\"oer. 
He constructed weak del Pezzo surfaces $X$ of characteristic two 
such that $H^1(X, \MO_X) \neq 0$ 
\cite[Theorem in Introduction]{Sch07}. 
Then Maddock discovered 
regular del Pezzo surfaces $X$ of characteristic two 
with $H^1(X, \MO_X) \neq 0$ \cite[Main Theorem]{Mad16}. 
If we allow singularities, 
it is known that 
there exists log del Pezzo surfaces $(X, \Delta)$ of characteristic three 
such that $H^1(X, \MO_X) \neq 0$ \cite{Tan1}. 

There are several results also in positive directions. 
Patakfalvi and Waldron proved that Gorenstein del Pezzo surfaces are geometrically normal when $p >3$ \cite[Theorem 1.5]{PW}. 
Fanelli and Schr\"oer showed that a regular del Pezzo surface $X$ 
is geometrically normal if $\rho(X)=1$ and the base field $k$ 
satisfies $[k:k^p] \leq 1$ \cite[Theorem 14.1]{FS}.  
Das proved that regular del Pezzo surfaces of characteristic $p\geq 5$ 
satisfy Kawamata--Viehweg vanishing \cite[Theorem 4.1]{Das}. 
Bernasconi and the author proved that log del Pezzo surfaces $(X, \Delta)$ of characteristic $p \geq 7$ 
are geometrically integral and satisfy $H^1(X, \MO_X)=0$ \cite[Theorem 1.7]{BT}.

\medskip

\textbf{Acknowledgements:} 
The author would like to thank Fabio Bernasconi and Gebhard Martin 
for useful comments. 
He also thanks the referee for many constructive suggestions and reading the manuscript carefully. 
The author was funded by 
the Grant-in-Aid for Scientific Research (KAKENHI No. 18K13386).

\section{Preliminaries}

\subsection{Notation}\label{ss-notation}

In this subsection, we summarise notation we will use in this paper. 

\begin{enumerate}
\item We will freely use the notation and terminology in \cite{Har77} 
and \cite{Kol13}. 
\item 
We say that a scheme $X$ is {\em regular} if 
the local ring $\MO_{X, x}$ at any point $x \in X$ is regular. 
\item 
For a scheme $X$, its {\em reduced structure} $X_{\red}$ 
is the reduced closed subscheme of $X$ such that the induced morphism 
$X_{\red} \to X$ is surjective. 
\item For an integral scheme $X$, 
we define the {\em function field} $K(X)$ of $X$ 
as $\MO_{X, \xi}$ for the generic point $\xi$ of $X$. 
\item 
For a field $k$, 
we say that $X$ is a {\em variety over} $k$ or a $k$-{\em variety} if 
$X$ is an integral scheme that is separated and of finite type over $k$. 
We say that $X$ is a {\em curve} over $k$ or a $k$-{\em curve} 
(resp. a {\em surface} over $k$ or a $k$-{\em surface}) 
if $X$ is a $k$-variety of dimension one (resp. two). 
\item 
For a variety $X$ over a field $k$, 
its normalisation is denoted by $X^N$. 
\item For a field $k$, we denote $\overline k$ an algebraic closure of $k$. 
If $k$ is of characteristic $p>0$, 
then we set $k^{1/p^{\infty}}:=\bigcup_{e=0}^{\infty} k^{1/p^e}
=\bigcup_{e=0}^{\infty} \{x \in \overline k\,|\, x^{p^e} \in k\}$. 
\item For an $\mathbb{F}_p$-scheme $X$ we denote by $F_X \colon X \to X$ the {\em absolute Frobenius morphism}. For a positive integer $e$ we denote by $F^e_X \colon X \to X$ 
the $e$-th iterated absolute Frobenius morphism.
\item If $k \subset k'$ is a field extension and $X$ is a $k$-scheme, we denote 
$X \times_{\Spec\,k} \Spec\,k'$ by $X \times_k k'$. 
\item 
Let $k$ be a field. 
A {\em del Pezzo surface} $X$ over $k$ 
is a projective normal surface over $k$ such that 
$-K_X$ is an ample $\Q$-Cartier divisor. 
\item 
Let $k$ be a field and let $X$ be a normal variety over $k$. 
We say that $X$ is {\em geometrically canonical} 
if $X \times_k \overline k$ is a normal variety over $\overline k$ 
which is canonical, i.e. has at worst canonical singularities. 
Note that if $X$ is geometrically canonical, then $X$ itself is canonical 
\cite[Proposition 2.3]{BT}. 
\item 
An $\F_p$-scheme $X$ is $F$-{\em finite} 
if the absolute Frobenius morphism $F:X \to X$ is a finite morphism. 
We say that a field $k$ of characteristic $p>0$ is $F$-{\em finite} 
if so is $\Spec\,k$, i.e. $[k:k^p]<\infty$. 
Note that if $k$ is an $F$-finite field and $X$ is of finite type over $k$, 
then also $X$ is $F$-finite. 
\item 
Let $X$ be a projective scheme over a field $k$ and 
let $F$ be a coherent sheaf on $X$. 
We say that $F$ is {\em globally generated} if 
there exist a positive integer $r$ and 
a surjective $\MO_X$-module homomorphism 
\[
\MO_X^{\oplus r} \to F. 
\]
An invertible sheaf $L$ on $X$ is {\em very ample over} $k$ if 
its complete linear system $|L|$ induces a closed immersion $X \hookrightarrow \mathbb P^N_k$. 
\end{enumerate}

\begin{dfn}[Definition {5.1} of \cite{Tan2}]\label{d-ell}
Let $k$ be a field of characteristic $p>0$ and 
let $X$ be a proper normal variety over $k$ with $H^0(X, \MO_X)=k$. 
Then we define the {\em Frobenius length of geometric non-normality} $\ell_F(X/k)$ 
of $X/k$ by 
{\small 
\[
\ell_F(X/k):=\min\{\ell \in \Z_{\geq 0}\,|\, (X \times_k k^{1/p^{\ell}})_{\red}^N
\text{ is geometrically normal over }k^{1/p^{\ell}}\},
\]
}
where the existence of the right hand side is guaranteed by 
\cite[Remark 5.2]{Tan2}. 
\end{dfn}

\begin{dfn}[Definition 7.4 of \cite{Tan2}]\label{d-epsilon}
Let $k$ be a field of characteristic $p>0$ and 
let $X$ be a proper normal variety over $k$ with $H^0(X, \MO_X)=k$. 
Set $R$ to be the local ring of $X \times_k k^{1/p^{\infty}}$ at the generic point. 
We define the {\em thickening exponent} $\epsilon(X/k)$ of $X/k$ by 
\[
\epsilon(X/k):=\log_p ({\rm length}_R R). 
\]
It follows from \cite[Theorem 7.3(1)]{Tan2} that $\epsilon(X/k)$ is a non-negative integer.  
\end{dfn}

\subsection{Summary of known results}

\begin{thm}\label{t-classify-bc}
Let $k$ be a field of characteristic $p>0$. 
Let $X$ be a canonical del Pezzo surface over $k$ 
such that $H^0(X, \MO_X)=k$. 
Set $Z:=(X \times_k \overline{k})_{\red}^N$. 
Then one of the following properties. 
\begin{enumerate}
\item $X$ is geometrically canonical over $k$. 
In particular, $Z=X \times_k \overline{k}$ and 
$Z$ is a canonical del Pezzo surface over $\overline k$. 
\item $X$ is not geometrically normal over $k$ and 
$Z \simeq \mathbb P_{\mathbb P^1}(\MO \oplus \MO(m))$ for some $m \in \Z_{\geq 0}$. 
\item $X$ is not geometrically normal over $k$ and 
$Z$ is isomorphic to a weighted projective surface $\mathbb P(1, 1, m)$ 
for some positive integer $m$. 
\end{enumerate} 
\end{thm}

\begin{proof}
See \cite[Theorem 3.3]{BT}. 
\end{proof}

\begin{thm}\label{t-ell-bdd}
Let $k$ be a field of characteristic $p>0$. 
Let $X$ be a canonical del Pezzo surface over $k$ 
such that $H^0(X, \MO_X)=k$. 
Then the following hold. 
\begin{enumerate}
\item 
If $p \geq 5$, then $X$ is geometrically canonical over $k$. 
\item 
If $p=3$, then $\ell_F(X/k)\leq 1$. 
\item 
If $p=2$, then $\ell_F(X/k) \leq 2$. 
\end{enumerate}
\end{thm}

\begin{proof}
See \cite[Theorem 3.7]{BT}. 
\end{proof}

\begin{thm}\label{t-T-PW}
Let $k$ be a field of characteristic $p>0$. 
Let $X$ be a proper normal variety over $k$ 
such that $H^0(X, \MO_X)=k$. 
Assume that $X$ is not geometrically normal over $k$. 
{Set $Z:=(X \times_k \overline{k})_{\red}^N$.} 
Then there exist nonzero effective $\Z$-divisors $C_1, ..., C_{\ell(X/k)}$ 
such that 
\[
K_Z + (p-1) \sum_{i=1}^{\ell(X/k)} C_i \sim f^*K_X
\]
where $f:Z \to X$ denotes the induced morphism. 
\end{thm}

\begin{proof}
See \cite[Proposition 5.11(2)]{Tan2}. 
\end{proof}

\begin{thm}\label{t-Frob-factor}
Let $k$ be a field of characteristic $p>0$. 
Let $X$ be a canonical del Pezzo surface over $k$ 
such that $H^0(X, \MO_X)=k$.  
Then the following hold. 
\begin{enumerate}
\item 
If $p=3$, then it holds that 
\[
(X \times_k k^{1/3})_{\red}^N \times_{k^{1/3}} \overline k 
\simeq (X \times_k \overline{k})_{\red}^N. 
\]
\item If $p=2$, then it holds that 
\[
(X \times_k k^{1/4})_{\red}^N \times_{k^{1/4}} \overline k 
\simeq (X \times_k \overline{k})_{\red}^N. 
\]
\end{enumerate}
\end{thm}

\begin{proof}
The assertion follows from Theorem \ref{t-ell-bdd} and \cite[Remark 5.2]{Tan2}. 
\end{proof}

\begin{thm}\label{t-Q-Gor}
Let $k$ be a field of characteristic $p>0$. 
Let $X$ be a regular del Pezzo surface over $k$
such that $H^0(X, \MO_X)=k$.  
Set $Z:=(X \times_k \overline k)_{\red}^N$. 
Then the following hold. 
\begin{enumerate}
\item 
If $p=3$, $3K_Z$ is Cartier. 
\item If $p=2$, then $4K_Z$ is Cartier. 
\end{enumerate}
\end{thm}

\begin{proof}
The assertion follows from Theorem \ref{t-ell-bdd} and \cite[Theorem 5.12]{Tan2}. 
\end{proof}

\section{Very ampleness}\label{s-va}

The purpose of this section is to prove that 
if $X$ is a regular del Pezzo surface, then 
$\omega_X^{-12}$ is very ample (Theorem \ref{t-dP-va2}).
To this end, we first establish a general criterion 
(Theorem \ref{t-va-criterion}) for very ampleness 
in Subsection \ref{ss1-va}. 
In Subsection \ref{ss2-va}, we apply this criterion to regular del Pezzo surfaces.

\subsection{A criterion for very ampleness}\label{ss1-va}

In this subsection, we give a criterion for very ampleness (Theorem \ref{t-va-criterion}). 
The strategy is a modification of 
Keeler's proof for base point freeness over algebraically closed fields \cite{Kee08}, 
which is in turn based on Smith's argument \cite{Smi97}. 
We first recall the definition (Definition \ref{d-0reg}) 
and a property (Lemma \ref{l-0reg-va}) of Castelnuovo--Mumford regularity.

\begin{dfn}\label{d-0reg}
Let $\kappa$ be an algebraically closed field. 
Let $Z$ be a projective scheme over $\kappa$. 
Let $H$ be an ample globally generated invertible sheaf on $Z$. 
A coherent sheaf $F$ on $Z$ is $0$-{\em regular} with respect to $H$ 
if 
\[
H^i(Z, F \otimes_{\MO_Z} H^{-i}) =0
\]
for any $i>0$. 
\end{dfn}

\begin{lem}\label{l-0reg-va}
Let $\kappa$ be an algebraically closed field. 
Let $Z$ be a projective scheme over $\kappa$ and let $z$ be a closed point on $Z$. 
Let $F$ be a coherent sheaf on $Z$ and let $H$ be an ample globally generated invertible sheaf on $Z$. 
Assume that $F$ is $0$-regular with respect to $H$. 
Then $F \otimes H \otimes \m_z$ is globally generated. 
\end{lem}

\begin{proof}
We may apply the same argument as in \cite[Lemma 3.7]{Wit17}. 
\end{proof}

\begin{thm}\label{t-va-criterion}
Fix a non-negative integer $d$. 
Let $k$ be a field of characteristic $p>0$. 
Let $X$ be a $d$-dimensional regular projective variety over $k$. 
Let $A$ be an ample invertible sheaf on $X$ and 
let $H$ be an ample globally generated invertible sheaf on $X$. 
Then $\omega_X \otimes_{\MO_X} H^{d+1} \otimes_{\MO_X} A$ is very ample over $k$. 
\end{thm}

\begin{proof}
We first reduce the problem to the case when $k$ is an $F$-finite field (cf. Subsection \ref{ss-notation}(12)). 
There exist a subfield $k_0 \subset k$, 
a projective scheme $X_0$ over $k_0$, 
and invertible sheaves $A_0$ and $H_0$ 
such that {$k_0$ is a field finitely generated over $\F_p$,} 
$X_0 \otimes_{k_0} k$, $f^*A_0=A$, and $f^*H_0=H$. 
Then we can check that {$k_0$ is $F$-finite and} $(k_0, X_0, A_0, H_0)$ 
satisfies the assumptions in the statement. 
Replacing $(k, X, A, H)$ by $(k_0, X_0, A_0, H_0)$, 
the problem is reduced to the case when $k$ is $F$-finite. 
In particular, also $X$ is $F$-finite (cf. Subsection \ref{ss-notation}(12)).

Fix $e \in \Z_{>0}$ and we denote the $e$-th iterated absolute Frobenius morphism 
$F^e:X \to X$ by $\Phi_e:X_e \to X$. 
Note that we consider $\Phi_e$ as a $k$-morphism, 
hence we distinguish $X$ and $X_e$ as $k$-schemes, 
although the equation $X_e=X$ holds as schemes. 
Let $A_e:=A$ and $H_e:=H$ be the invertible sheaves on $X_e$. 
Note that we have $\Phi_e^*A=A_e^{p^e}$ and $\Phi_e^*H=H_e^{p^e}$. 

For $\kappa:=\overline k$, we take the base changes
\[
\begin{CD}
Y_e @>\alpha_e >> X_e\\
@VV\Psi_e V @VV\Phi_e V\\
Y @>\alpha >> X\\
@VVV @VVV\\
\Spec\,\kappa @>>> \Spec\,k,
\end{CD}
\]
hence both the above squares are cartesian. 
We set $A':=\alpha^*A, H':=\alpha^*H$, $A'_e:=\alpha_e^*A_e$, and $H'_e := \alpha_e^*H_e$. 
Since $\Phi_e^*A=A_e^{p^e}$ and $\Phi_e^*H=H_e^{p^e}$, 
we have $\Psi_e^*A'=A_e'^{p^e}$ and $\Psi_e^*H'=H_e'^{p^e}$.

\begin{claim}\label{c-va-criterion}
There exists a positive integer $e$ such that 
the coherent sheaf $(\Psi_e)_*(\omega_{Y_e} \otimes H_e'^{p^ed} \otimes A_e'^{p^e})$ 
on $Y$ is $0$-regular with respect to $H'$, i.e. the equation 
\[
H^i(Y, (\Psi_e)_*(\omega_{Y_e} \otimes H_e'^{p^ed} \otimes A_e'^{p^e})  \otimes H'^{-i})=0
\]
holds for any $i>0$. 
\end{claim}

\begin{proof}(of Claim \ref{c-va-criterion})
We have 
\begin{eqnarray*}
&&
H^i(Y, (\Psi_e)_*(\omega_{Y_e} \otimes H_e'^{p^ed} \otimes A_e'^{p^e})  \otimes H'^{-i})\\
& \simeq &
H^i(Y, (\Psi_e)_*(\omega_{Y_e} \otimes H_e'^{p^e(d-i)} \otimes A_e'^{p^e}))\\
& \simeq &
H^i(Y_e, \omega_{Y_e} \otimes H_e'^{p^e(d-i)} \otimes A_e'^{p^e}),
\end{eqnarray*}
where the first isomorphism follows from the projection formula and 
the second isomorphism holds because $\Psi_e$ is an affine morphism. 
By flat base change theorem, it holds that 
\[
H^i(Y_e, \omega_{Y_e} \otimes H_e'^{p^e(d-i)} \otimes A_e'^{p^e})
\simeq 
H^i(X_e, \omega_{X_e} \otimes H_e^{p^e(d-i)} \otimes A_e^{p^e}) \otimes_k \kappa. 
\]
Recall that $X$ and $X_e$ are isomorphic as schemes.  
Therefore, we have an isomorphism as abelian groups: 
\[
H^i(X_e, \omega_{X_e} \otimes H_e^{p^e(d-i)} \otimes A_e^{p^e})
\simeq 
H^i(X, \omega_X \otimes (H^{d-i} \otimes A)^{p^e}). 
\]
It is enough to treat the case when $i \leq  \dim X = d$. 
Hence, $H^{d-i} \otimes A$ is ample. 
Then, by the Serre vanishing theorem, 
the right hand side is equal to zero for $e \gg 0$. 
This completes the proof of Claim \ref{c-va-criterion}. 
\end{proof}

Fix a closed point $y$ on $Y$. 
Take a positive integer $e$ as in Claim \ref{c-va-criterion}. 
Then $(\Psi_e)_*(\omega_{Y_e} \otimes H_e'^{p^ed} \otimes A_e'^{p^e})$ 
is $0$-regular with respect to $H'$. 
Lemma \ref{l-0reg-va} implies that the coherent sheaf 
\[
(\Psi_e)_*(\omega_{Y_e}) \otimes H'^{d+1} \otimes A' \otimes \m_y
=
(\Psi_e)_*(\omega_{Y_e} \otimes H_e'^{p^ed} \otimes A_e'^{p^e}) \otimes H' \otimes \m_y
\] 
is globally generated. 

Since $\MO_{X, x} \to ((\Phi_e)_*\MO_{X_e})_x$ splits for any point $x$ on $X$ 
\cite[Theorem 107 in Section 42]{Mat80}, 
we obtain a surjective $\MO_X$-module homomorphism 
$(\Phi_e)_*(\omega_{X_e}) \to \omega_X$ by applying $\mathcal Hom_{\MO_X}(-, \omega_X)$ 
to $\MO_X \to (\Phi_e)_*\MO_{X_e}$. 
Taking the base change $(-) \times_k \kappa$, 
there exists a surjective $\MO_Y$-module homomorphism 
$(\Psi_e)_*(\omega_{Y_e}) \to \omega_Y$, 
which induces another surjective $\MO_Y$-module homomorphism 
\[
(\Psi_e)_*(\omega_{Y_e}) \otimes H'^{d+1} \otimes A' \otimes \m_y
\to 
\omega_Y \otimes H'^{d+1} \otimes A' \otimes \m_y.
\]
Since $(\Psi_e)_*(\omega_{Y_e}) \otimes H'^{d+1} \otimes A' \otimes \m_y$ 
is globally generated, 
also $\omega_Y \otimes H'^{d+1} \otimes A' \otimes \m_y$ 
is globally generated. 
This implies that $\omega_Y \otimes H'^{d+1} \otimes A'$ is very ample over $\kappa$. 
Since very ampleness descends by base changes, 
$\omega_X \otimes H^{d+1} \otimes A$ is very ample over $k$. 

{In what follows, we prove that the very ampleness acutually descends. 
First of all, the base point freeness descends, because it is characterised by the surjectivity of $H^0(X, L) \otimes_k \MO_X \to L$ for $L := \omega_X \otimes H^{d+1} \otimes A$. 
We then get the  morphism $\varphi : X \to \mathbb P^N_k$ induced by $|L|$, where $N := \dim_k H^0(X, L)-1$. 
Since the base change $\varphi \times_k \kappa$ is a closed immersion, $\varphi$ is a finite morphism, because there exists no curve contracted by $\varphi$. 
Hence it is enough to show that, given a ring homomorphism $\psi : A \to B$ of $k$-algebra, 
$\psi$ is surjective if so is $\psi \otimes_k \kappa$. This follows from the fact that  
$k \hookrightarrow \kappa$ is faithfully flat.}
\end{proof}

\subsection{Very ampleness for regular del Pezzo surfaces}\label{ss2-va}

In this subsection, we prove the main result (Theorem \ref{t-dP-va2}) 
of this section. 
We first focus on the case when $X$ is not geometrically normal. 

\begin{thm}\label{t-dP-va1}
Let $k$ be a field of characteristic $p>0$. 
Let $X$ be a regular del Pezzo surface over $k$ such that $H^0(X, \MO_X)=k$. 
Let $A$ be an ample invertible sheaf and 
let $N$ be a nef invertible sheaf. 
Assume that $X$ is not geometrically normal over $k$. 
Then the following hold. 
\begin{enumerate}
\item If $p=2$, then $A^4$ is globally generated. 
\item If $p=3$, then $A^3$ is globally generated. 
\item If $p=2$, then $\omega_X^{-12} \otimes N$ is very ample over $k$. 
\item If $p=3$, then $\omega_X^{-9} \otimes N$ is very ample over $k$. 
\end{enumerate}
\end{thm}

\begin{proof}
If $p=2$, then we set $e:=2$ and $q:=p^e=4$. 
If $p=3$, then we set $e:=1$ and $q:=p^e=3$.

Let us prove that $A^q$ is globally generated. 
Set $A_{\overline k}$ to be the pullback of $A$ to $X \times_k \overline k$. 
Since $e \geq \ell_F(X/k)$, the $e$-th iterated absolute Frobenius factors 
(Theorem \ref{t-Frob-factor}): 
\[
F^e_{X \times_k \overline k}:X \times_k \overline k \xrightarrow{\psi} Z=(X \times_k \overline k)_{\red}^N 
\xrightarrow{\varphi} X \times_k \overline k. 
\]
Thus, $\varphi^*(A_{\overline k})$ is an ample invertible sheaf on a projective toric surface $Z$ (Theorem \ref{t-classify-bc}). 
Then $\varphi^*(A_{\overline k})$ is globally generated, 
hence so is its pullback:  
\[
\psi^*\varphi^*(A_{\overline k})
=(F^e_{X \times_k \overline k})^*(A_{\overline k})
=A_{\overline k}^q.
\]
Hence, also $A_{\overline k}^q$ is globally generated. 
Thus, (1) and (2) hold. 

Let us prove (3) and (4). 
By (1) and (2), $\omega_X^{-q}$ is globally generated. 
Then it follow from Theorem \ref{t-va-criterion} 
that the invertible sheaf 
\[
\omega_X^{-3q} \otimes N = \omega_X \otimes (\omega_X^{-q})^{\dim X+1} \otimes 
(\omega_X^{-1} \otimes N) 
\] 
is very ample over $k$. 
Thus (3) and (4) hold. 
\end{proof}

\begin{thm}\label{t-dP-va2}
Let $k$ be a field of characteristic $p>0$. 
Let $X$ be a regular del Pezzo surface over $k$ such that $H^0(X, \MO_X)=k$. 
Then $\omega_X^{-m}$ is very ample over $k$ for any integer $m$ such that $m \geq 12$. 
\end{thm}

\begin{proof}
If $X$ is not geometrically normal over $k$, then 
the assertion follows from Theorem \ref{t-dP-va1}. 
Assume that $X$ is geometrically normal over $k$. 
Then $X$ is geometrically canonical over $k$ (Theorem \ref{t-classify-bc}). 
In this case, $\omega_X^{-2}$ is globally generated by \cite[Proposition 2.14(1)]{BT}. 
Hence, it follows from Theorem \ref{t-va-criterion} that $\omega_X^{-m} = \omega_X \otimes (\omega_X^{-2})^3 \otimes \omega_X^{-(m-5)}$ is very ample for $m \geq 6$. 
\end{proof}

\section{Boundedness of volumes}\label{s-bdd-vol}

The purpose of this section is to 
show Theorem \ref{t-vol-bound2}, which gives the inequality 
\[
K_X^2 \leq \max\{9, 2^{2r+1}\}
\]
for a regular del Pezzo surface $X$ over a field $k$ of characteristic $p>0$ 
such that $H^0(X, \MO_X)=k$ and $r:=\log_p [k:k^p]$. 
If $X$ is geometrically normal, then the problem has been settled already 
(cf. the proof of Theorem \ref{t-bd-not-gn2}(1)). 
Most of this subsection is devoted to 
analysis of the geometrically non-normal case. 
In Subsection \ref{ss1-bdd-vol}, 
we first restrict possibilities for $(Z, D)$, where 
$Z:=(X \times_k \overline k)_{\red}^N$ and $D$ is an effective divisor $D$ on $Z$ such that the linear equivalence 
\[
K_Z+D \sim g^*K_X 
\]
holds for the induced morphism $g:Z \to X$. 
In Subsection \ref{ss2-bdd-vol}, 
we prove that there are only finitely many possibilities for $K_X^2$ 
after we fix $\epsilon(X/k)$ (Theorem \ref{t-bd-not-gn}). 
We then obtain our main result (Theorem \ref{t-vol-bound2}) 
by combining with fundamental properties on $\epsilon(X/k)$.

\subsection{Restriction on possibilities}\label{ss1-bdd-vol}

The purpose of this subsection is to prove the following proposition.

\begin{prop}\label{p-restriction}
Let $k$ be a field of characteristic $p>0$. 
Let $X$ be a canonical del Pezzo surface over $k$ such that $H^0(X, \MO_X)=k$. 
Set $Z:=(X \times_k \overline k)_{\red}^N$ 
and let $g:Z \to X$ be the induced morphism. 
If $D$ is a nonzero effective divisor on $Z$ satisfying
\begin{equation}\label{e1-restriction}
K_Z + D \sim g^*K_X,
\end{equation} 
then one of the following {holds}. 
\begin{enumerate}
\item $Z \simeq \mathbb P^2$. In this case, it holds that 
\begin{enumerate}
\item $\MO_Z(D) \simeq \MO(1)$, or 
\item $\MO_Z(D) \simeq \MO(2)$
\end{enumerate}
\item $Z \simeq \mathbb P^1 \times \mathbb P^1$. In this case, it holds that
\begin{enumerate}
\item $\MO_Z(D) \simeq \MO(1, 1)$,  
\item $\MO_Z(D) \simeq \MO(1, 0)$, or
\item $\MO_Z(D) \simeq \MO(0, 1)$. 
\end{enumerate}
\item 
$Z \simeq \mathbb P(1, 1, m)$ for some $m \geq 2$. 
In this case, $D \sim 2F$, where $F$ is a prime divisor such that 
$F^2 =1/m $.  
\item $Z \simeq \mathbb P( \MO \oplus \MO(m))$ for some $m \geq 1$. 
In this case, if $\pi:Z \to \mathbb P^1$ is the $\mathbb P^1$-bundle structure, 
$F$ is a fibre of $\pi$,  and $C$ is a curve with $C^2=-m$, then  
\begin{enumerate}
\item $D \sim C$, or 
\item $D \sim C+F$. 
\end{enumerate}
\end{enumerate}
\end{prop}

\begin{proof}
Note that $-(K_Z+D)$ is ample. 
Hence, if $Z \simeq \mathbb P^2$ or $Z \simeq \mathbb P^1 \times \mathbb P^1$, 
then (1) or (2) holds. 
We assume that $Z$ is isomorphic to neither $\mathbb P^2$ nor $\mathbb P^1 \times \mathbb P^1$. 
Then it follows from Theorem \ref{t-classify-bc} that 
there is $m \geq 1$ such that either 
\begin{enumerate}
\item[(i)] $Z \simeq \mathbb P_{\mathbb P^1}(\MO \oplus \MO(m))$, or 
\item[(ii)] $Z \simeq \mathbb P(1, 1, m)$ and $m \geq 2$. 
\end{enumerate}

Then, for the minimal resolution $\mu:W \to Z$,  
it holds that $W \simeq \mathbb P_{\mathbb P^1}(\MO \oplus \MO(m))$. 
We have the induced morphisms: 
\[
h:W \xrightarrow{\mu} Z \xrightarrow{g} X
\]
Let $\pi:W \to \mathbb P^1$ be the $\mathbb P^1$-bundle structure. 
Let $F_W$ be a fibre of $\pi$ and let $C$ be the curve on $W$ such that $C^2=-m$. 
For $D_W:=\mu_*^{-1}D$, we obtain 
\begin{equation}\label{e2-restriction}
K_W+D_W+c C \sim \mu^*(K_Z+D) \sim h^*K_X
\end{equation}
for some $c \in \Z_{\geq 0}$. 
We have 
\[
-K_W \sim 2C+(m+2)F_W 
\]
and $D_W \sim aC+bF_W$ for some $a, b \in \Z_{ \geq 0}$ with $(a, b) \neq (0, 0)$. 
Thus it holds that 
\begin{equation}\label{e3-restriction}
-h^*K_X \sim -K_W-D_W-c C \sim (2-a-c)C + (m+2-b)F_W.
\end{equation}

We first show that $a+c=1$. 
Since $-h^*K_X$ is big, we obtain $(-h^*K_X) \cdot F_W>0$, hence it holds that $2-a-c \geq 1$. 
Then we have $1 \leq 2-a-c \leq 2$. 
Thus, it is enough to prove that $a+c \neq 0$. 
Assuming $a=c=0$, let us derive a contradiction. 
We have 
\[
-h^*K_X \cdot C= (2C + (m+2-b)F_W) \cdot C =-m+2-b
\]
If (i) holds, then $-h^*K_X$ is ample, hence we obtain 
$0 < -h^*K_X \cdot C =-m+2-b \leq 1-b$, which  in turn implies $b=0$. 
If (ii) holds, then it holds that  
$0 \leq  -h^*K_X \cdot C =-m+2-b \leq -b$. 
In any case, we have $b=0$, which contradicts $(a, b) \neq (0, 0)$. 
Therefore, we obtain $a+c=1$. 
In particular, (\ref{e3-restriction}) implies that 
\begin{equation}\label{e4-restriction}
-h^*K_X \sim C + (m+2-b)F_W.
\end{equation}
We treat the following two cases separately: 
\[
(a, c)=(0, 1) \quad \text{or}\quad (a, c)=(1, 0). 
\]

Let us handle the case when $(a, c)=(0, 1)$. 
By $c \neq 0$, $\mu$ is not an isomorphism, 
hence we obtain $Z \simeq \mathbb P(1, 1, m)$. 
Since $h^*K_X \cdot C=\mu^*(K_Z+D)\cdot C=0$, (\ref{e4-restriction}) implies $b=2$. 
Thus, we conclude $(a, b, c)=(0, 2, 1)$. 
This implies that (3) holds.

Then we may assume that $(a, c)=(1, 0)$. 
Assume (i). 
Then $-h^*K_X$ is ample. 
By (\ref{e4-restriction}), we have 
\[
0< -h^*K_X \cdot C= (C + (m+2-b)F_W) \cdot C=-m+(m+2-b)=2-b. 
\]
Therefore, we obtain $b \in \{0, 1\}$. Thus, (4) holds. 
Assume (ii). 
Since $c=0$ and $c$ is defined by (\ref{e2-restriction}), 
we have $m=2$. 
Again by (\ref{e2-restriction}), we obtain 
\[
0= h^*K_X \cdot C = (K_W+D_W+c C) \cdot C = D_W \cdot C = (aC+bF_W) \cdot C = -2+b.
\]
Thus, it holds that $(a, b, c)=(1, 2, 0)$. Thus, (3) holds. 
\end{proof}

\begin{rem}\label{r-restriction}
We use notation as in Proposition \ref{p-restriction}. 
Note that $(g^*K_X)^2=(h^*K_X)^2$. 
By direct computation using (\ref{e4-restriction}), the following hold. 
\begin{enumerate}
\item If $Z \simeq \mathbb P^2$, then $(g^*K_X)^2 \in \{1, 4\}$. 
\item If $Z \simeq \mathbb P^1 \times \mathbb P^1$, then $(g^*K_X)^2 \in \{2, 4\}$. 
\item If $Z \simeq \mathbb P(1, 1, m)$ for some $m \geq 2$, then $(g^*K_X)^2=m$. 
\item If $Z \simeq \mathbb P_{\mathbb P^1}( \MO \oplus \MO(m))$ for some $m \geq 1$, 
then $(g^*K_X)^2 \in \{m+2, m+4\}$. 
\end{enumerate}
\end{rem}

\begin{rem}\label{r-restriction2}
We use notation as in Proposition \ref{p-restriction}. 
If $p=3$, then we can find a nonzero effective divisor $D'$ such that 
$K_Z+2D' \sim g^*K_X$. 
In this case, (2) and (4) in Proposition \ref{p-restriction} does not occur. 
\end{rem}

\subsection{Classification of base changes}\label{ss2-bdd-vol}

In this subsection \ref{t-vol-bound2}, 
we prove the main result of this section (Theorem \ref{t-vol-bound2}), 
which asserts the inequality 
\[
K_X^2 \leq \max\{9, 2^{2r+1}\}
\]
for a regular del Pezzo surface $X$ over a field $k$ of characteristic $p>0$ 
such that $H^0(X, \MO_X)=k$ and $r:=\log_p [k:k^p]$. 
This result is a consequence of the boundedness of $K_X^2$ 
in terms of $\epsilon(X/k)$ (Theorem \ref{t-bd-not-gn}). 
To this end, we prove a kind of classification 
after the base change to the algebraic closure (Theorem \ref{t-bd-not-gn}). 
We first establish auxiliary results: Lemma \ref{l-pic-one} 
and Lemma \ref{l-epsilon-bc}.

\begin{lem}\label{l-pic-one}
Let $X$ be a regular del Pezzo surface over $k$ such that $H^0(X, \MO_X)=k$. 
Set $Z:=(X \times_k \overline k)_{\red}^N$. 
Assume that $\rho(X)=1$ and $X$ is not geometrically normal. 
Then it holds that  $Z \simeq \mathbb P^1 \times \mathbb P^1$ 
or $Z \simeq \mathbb P(1, 1, m)$ for some $m \geq 1$.  
\end{lem}

\begin{proof}
Assume that $Z$ is not isomorphic to $\mathbb P(1, 1, m)$, then 
it follows from Theorem \ref{t-classify-bc} 
that $Z \simeq \mathbb P_{\mathbb P^1}( \MO \oplus \MO(m))$ 
for some $m \geq 0$. 
Suppose $m>0$ and let us derive a contradiction. 
Set $\kappa:= k^{1/p^{\infty}}$ and $Y:=(X \times_k \kappa)_{\red}^N$. 
Then we have $Y \times_{\kappa} \overline k \simeq Z$. 
Hence, $Y$ is smooth over $\kappa$. 
We have $\rho(Y)=1$ \cite[Proposition 2.4(3)]{Tan18b}.

Let $\pi:Z \to B$ be the $\mathbb P^1$-bundle structure. 
There is a finite Galois extension $\kappa'/\kappa$ 
such that $\pi$ descends to $\kappa'$, 
i.e. there exists a $\kappa'$-morphism $\pi':Z' \to B'$ of 
smooth $\kappa'$-varieties whose base change by $(-) \times_{\kappa'} \overline k$ is $\pi:Z \to B$. 
Let $F'$ be a fibre of $\pi'$ over a closed point. 
For the Galois group $G$ of $\kappa'/\kappa$ and any element $\sigma \in G$, 
we have that $\sigma^*(F')^2=F'^2=0$. 
If $\sigma^*(F')$ is not a fibre of $\pi'$, 
then $\sigma^*(F)$ induces another fibration which deduces that 
$Z \simeq \mathbb P^1 \times \mathbb P^1$. 
Hence $\sigma^*(F)$ is a fibre of $\pi'$. 
Then $\widetilde{F}:=\sum_{\sigma \in G} \sigma^*(F)$ satisfies $\widetilde F^2=0$. 
As $\widetilde F$ descends to $Y$, 
there exists an effective divisor $D$ on $Y$ such that $D^2=0$. 
However, this contradicts $\rho(Y)=1$. 
\end{proof}

\begin{lem}\label{l-epsilon-bc}
Let $X$ be a regular del Pezzo surface over $k$ such that $H^0(X, \MO_X)=k$. 
Set $Z:=(X \times_k \overline k)_{\red}^N$ 
and 
let $g:Z = ((X \times_k \overline k)_{\red}^N \to X$ be the induced morphism. 
Then it holds that $p^{\epsilon(X/k)}(g^*K_X)^2=K_X^2$.  
\end{lem}

\begin{proof}
The assertion follows from 
Definition \ref{d-epsilon} and 
\cite[Example 1 in page 299]{Kle66}.
\end{proof}

\begin{thm}\label{t-bd-not-gn}
Let $X$ be a regular del Pezzo surface over $k$ such that $H^0(X, \MO_X)=k$. 
Set $Z:=(X \times_k \overline k)_{\red}^N$ and let $g:Z \to X$ be the induced morphism. 
Assume that $X$ is not geometrically normal over $k$. 
Then there exists a nonzero effective $\Z$-divisor $E$ on $Z$ such that 
\begin{equation}\label{e1-bd-not-gn}
K_Z + (p-1)E \sim g^*K_X.
\end{equation} 
Furthermore, if $E$ is a nonzero effective divisor $E$ on $Z$ 
satisfying (\ref{e1-bd-not-gn}), then the following hold. 
\begin{enumerate}
\item It holds that $p=2$ or $p=3$.  
\item 
If $p=3$, then the quadruple $(Z, E, (g^*K_X)^2, K_X^2)$ 
satisfies one of the possibilities in the following table. 
\begin{table}[H]
\caption{$p=3$ case}
\begin{tabular}{|c|c|c|c|} \hline
$Z$ & $E$ & $(g^*K_X)^2$ & $K_X^2$ \\ \hline \hline
$\mathbb P^2$  & $\MO(1)$ & $1$ & $3^{\epsilon(X/k)}$   \\ \cline{2-4} 
\hline 
$\mathbb P(1, 1, 3)$ & $F$ & $3$ &$3^{\epsilon(X/k)+1}$ \\ \cline{2-4}
\hline
  \end{tabular}
\end{table}
\item 
If $p=2$, then the quadruple $(Z, D, (g^*K_X)^2, K_X^2)$ 
satisfies one of the possibilities in the following table. 
\begin{table}[H]
\caption{$p=2$ case}
\begin{tabular}{|c|c|c|c|} \hline
$Z$ & $E$ & $(g^*K_X)^2$ & $K_X^2$ \\ \hline \hline
$\mathbb P^2$ & $\MO(1)$ & $4$ & $2^{\epsilon(X/k)+2}$  \\ \cline{2-4}
      & $\MO(2)$ & $1$ & $2^{\epsilon(X/k)}$  \\ 
\hline
$\mathbb P(1, 1, 2)$ & $2F$ & $2$ &$2^{\epsilon(X/k)+1}$ \\ \cline{2-4}
\hline
$\mathbb P(1, 1, 4)$ & $2F$ & $4$ &$2^{\epsilon(X/k)+2}$ \\ \cline{2-4}
\hline 
$\mathbb P^1 \times \mathbb P^1$     & $\MO(1, 0)$ & $4$ & $2^{\epsilon(X/k)+2}$ \\ 
   & $\MO(1, 1)$ & $2$ &$2^{\epsilon(X/k)+1}$ \\ \cline{2-4}
\hline
$\mathbb P_{\mathbb P^1}(\MO \oplus \MO(1))$       & $C$ & $5$ &$5 \cdot 2^{\epsilon(X/k)}$ \\ \cline{2-4}
 & $C+F$ & $3$ & $3 \cdot 2^{\epsilon(X/k)}$ \\ \cline{2-4}
\hline
$\mathbb P_{\mathbb P^1}(\MO \oplus \MO(2))$       & $C$ & $6$ &$3 \cdot 2^{\epsilon(X/k)+1}$ \\ \cline{2-4}
 & $C+F$ & $4$ & $2^{\epsilon(X/k)+2}$ \\ \cline{2-4}
\hline
$\mathbb P_{\mathbb P^1}(\MO \oplus \MO(4))$       & $C$ & $8$ &$2^{\epsilon(X/k)+3}$ \\ \cline{2-4}
 & $C+F$ & $6$ & $3 \cdot 2^{\epsilon(X/k)+1}$ \\ \cline{2-4}
\hline
\end{tabular}
\end{table}
\end{enumerate}
Here, if we write an invertible sheaf in the list, then 
it means that $\MO_Z(E)$ is isomorphic to it. 
If we write a divisor, then it means that $E$ is linearly equivalent to it. 
On $\mathbb P(1, 1, m)$ with $m \geq 2$, $F$ denotes a prime divisor such that $F^2=1/m$. 
On $\mathbb P_{\mathbb P^1}(\MO \oplus \MO(m))$ with $m \geq 1$, 
$C$ is the curve such that $C^2=-m$ and 
$F$ denotes a fibre of the $\mathbb P^1$-bundle structure 
$\mathbb P_{\mathbb P^1}(\MO \oplus \MO(m)) \to \mathbb P^1$. 
\end{thm}

\begin{proof}
The existence of $E$ follows from Theorem \ref{t-T-PW}.  
The assertion (1) holds by \cite[Theorem 1.5]{PW}. 
We omit the proof of (2), as 
it is similar and easier than the one of (3). 

Let us show (3). 
Pick a nonzero effective divisor $E$ on $Z$ satisfying (\ref{e1-bd-not-gn}). 
If $(Z, E)$ is one of the possibilities in the table, 
then $(g^*K_X)^2$ and $K_X^2$ automatically determined. 
Thus, it is enough to show that the pair $(Z, E)$ satisfies one of the possibilities.

We first treat the following two cases: 
\begin{enumerate}
\item[(i)] $Z \simeq \mathbb P^1 \times \mathbb P^1$. 
\item[(ii)] $Z \simeq \mathbb P(1, 1, m)$ for some $m \geq 1$. 
\end{enumerate}
If (i) holds, then Remark \ref{r-restriction}(2) implies the assertion. 
Assume that (ii) holds. 
If $m=1$, then the assertions follow from Remark \ref{r-restriction}(1)
Let us handle the case when $m \geq 2$. 
It follows from Theorem \ref{t-Q-Gor} that $m$ is a divisor of  $4$. 
By Remark \ref{r-restriction}(3), the assertion holds. 
In particular, by Lemma \ref{l-pic-one}, 
we are done for the case when $\rho(X)=1$. 

We now treat the case when $\rho(X) \neq 1$. 
We have $\rho(X) \leq \rho(Z) \leq 2$, 
where the latter inequality follows from Proposition \ref{p-restriction}. 
Hence we have $\rho(X)=\rho(Z)=2$. 
Since the case (i) has been settled already, 
Proposition \ref{p-restriction} enables us to assume that 
the case (4) of Proposition \ref{p-restriction} occurs, 
i.e. $Z \simeq \mathbb P_{\mathbb P^1}(\MO \oplus \MO(m))$ for some $m \geq 1$. 
By \cite[Theorem 4.4]{Tan18a}, 
there are two extremal contractions $\varphi:X \to X'$ and $X \to X''$. 
Both of them induce morphisms $Z \to Z'$ and $Z \to Z''$ 
with $\dim X'=\dim Z'$ and $\dim X''=\dim Z''$. 
Hence we may assume that $\dim X'=2$, i.e. 
$\varphi:X \to X'$ is a birational morphism that contracts a single curve. 
Then $X'$ is a regular del Pezzo surface with $\rho(X')=1$ {\cite[Theorem 10.5]{Kol13}}. 

Assume that $X'$ is not geometrically normal. 
Then $Z' \simeq \mathbb P^1 \times \mathbb P^1$ or $Z' \simeq \mathbb P(1, 1, m)$ 
for some $m \in \{1, 2, 4\}$ (Lemma \ref{l-pic-one}). 
Hence, we may assume that $X'$ is geometrically normal. 
Then $X'$ is geometrically canonical (Theorem \ref{t-classify-bc}). 
Therefore, $Z'$ has at worst canonical singularities.  
In particular, we obtain $m \leq 2$. 
Hence, Remark \ref{r-restriction}(4) implies the assertion. 
\end{proof}

\begin{thm}\label{t-bd-not-gn2}
Let $k$ be a field of characteristic $p>0$. 
Let $X$ be a regular del Pezzo surface over $k$ such that $H^0(X, \MO_X)=k$. 
Then the following hold. 
\begin{enumerate}
\item If $X$ is geometrically normal, then $K_X^2 \leq 9$. 
\item 
Assume that 
$X$ is not geometrically normal. 
Then $p \in \{2, 3\}$ and the following hold. 
\begin{enumerate}
\item If $p =3$, then $K_X^2 \leq 3^{\epsilon(X/k)+1}$. 
\item If $p=2$, then $K_X^2 \leq 2^{\epsilon(X/k)+3}$. 
\end{enumerate}
\end{enumerate}
In particular, if $X$ is geometrically reduced, then 
it holds that $K_X^2 \leq 9$. 
\end{thm}

\begin{proof}
Let us show (1). 
If $X$ is geometrically normal, then $X$ is geometrically canonical (Theorem \ref{t-classify-bc}). 
Hence, we have $K_X^2 \leq 9$ (cf. \cite[Lemma 5.1]{BT}). 
Thus (1) holds. 

Let us show (2). 
Assume that 
$X$ is not geometrically normal. 
Then \cite{PW} implies that $p \in \{2, 3\}$. 
The assertions (a) and (b) follow directly from Theorem \ref{t-bd-not-gn}. 
Note that the last assertion holds by the fact that 
$\epsilon(X/k)=0$ if $X$ is geometrically reduced over $k$ (Definition \ref{d-epsilon}). 
\end{proof}

\begin{cor}\label{c-bd-not-gn2}
Let $k$ be a field of characteristic $p>0$. 
Let $X$ be a regular del Pezzo surface over $k$ such that $H^0(X, \MO_X)=k$. 
Then the following hold. 
\begin{enumerate}
\item If $p \geq 5$, then $K_X^2 \leq 9$. 
\item If $p =3$, then $K_X^2 \leq \max \{9, 3^{\epsilon(X/k)+1}\}$. 
\item If $p =2$, then $K_X^2 \leq \max \{9, 2^{\epsilon(X/k)+3}\}$. 
\end{enumerate}
\end{cor}

\begin{proof}
The assertion follows from Theorem \ref{t-bd-not-gn2}. 
\end{proof}

\begin{thm}\label{t-vol-bound2}
Let $k$ be a field of characteristic $p>0$ such that $[k:k^p]<\infty$. 
Let $X$ be a regular del Pezzo surface over $k$ such that $H^0(X, \MO_X)=k$. 
Then the following hold. 
\begin{enumerate}
\item If $p \geq 5$, then $K_X^2 \leq 9$. 
\item If $p=3$, then $K_X^2 \leq \max\{9, [k:k^3]\}$. 
\item If $p=2$, then $K_X^2 \leq \max\{9, 2\cdot ([k:k^2])^2\}$. 
\end{enumerate}
In particular, if $r:=\log_p[k:k^p]$, then it holds that 
\[
K_X^2 \leq \max \{9, 2^{2r+1}\}. 
\]
\end{thm}

\begin{proof}
If $X$ is geometrically normal, then $K_X^2 \leq 9$ (Theorem \ref{t-bd-not-gn2}). 
Hence we may assume that $X$ is not geometrically normal. 
In this case, we have $p \in \{2, 3\}$ (Theorem \ref{t-bd-not-gn2}). 
Hence, (1) holds. 

Since $X$ is not geometrically normal, we have $[k:k^p] \neq 1$. Hence, 
it follows from \cite[Remark 1.7]{Tan2} that 
\[
\epsilon(X/k) \leq \ell_F(X/k) (\log_p[k:k^p]-1). 
\]
In particular, we have that 
\[
p^{\epsilon(X/k)} \leq p^{\ell_F(X/k) (\log_p[k:k^p]-1)}=(p^{-1} \cdot [k:k^p])^{\ell_F(X/k)}. 
\]

Let us show (2). 
We have $\ell_F(X/k) \leq 1$ (Theorem \ref{t-ell-bdd}) 
and $K_X^2 \leq 3^{\epsilon(X/k)+1}$ (Theorem \ref{t-bd-not-gn2}). 
Therefore, we obtain 
\[
K_X^2 \leq 3^{\epsilon(X/k)+1} \leq 3 \cdot (3^{-1} \cdot [k:k^3])^{\ell_F(X/k)} 
\leq [k:k^3]. 
\]
Thus (2) holds.

Let us show (3). 
We have $\ell_F(X/k) \leq 2$ (Theorem \ref{t-ell-bdd}) 
and $K_X^2 \leq 2^{\epsilon(X/k)+3}$ (Theorem \ref{t-bd-not-gn2}). 
Therefore, we obtain 
\[
K_X^2 \leq 2^{\epsilon(X/k)+3} \leq 2^3 \cdot (2^{-1} \cdot [k:k^2])^{\ell_F(X/k)} 
\leq 2 \cdot ([k:k^2])^2. 
\]
Thus (3) holds. 
\end{proof}

\section{Boundedness of regular del Pezzo surfaces}\label{s-bounded}

In this section, 
we prove the boundedness of geometrically integral regular del Pezzo surfaces (Theorem \ref{t-dP-bdd3}). 
The proof will be given in Subsection \ref{ss2-bounded}. 
In Subsection \ref{ss1-bounded}, we recall results on Chow varieties.



\subsection{Chow varieties}\label{ss1-bounded}

The purpose of this subsection 
is to give a proof of Proposition \ref{p-chow-var}. 
The result itself is well known to experts, 
however we give a proof for the sake of completeness. 
Since we shall use Chow varieties, 
we now recall its construction and results for later use
\cite[Ch. I, Section 3, Section 4]{Kol96}.

\begin{dfn}
Let $Chow_{r, d}(\mathbb P^N/\Z)$ be the contravariant functor 
from the category of semi-normal schemes to the category of sets such that 
if $T$ is a semi-normal scheme, then 
$Chow_{r, d}(\mathbb P^N/\Z)(T)$ is the set 
of well-defined algebraic families of nonnegative cycles of $\mathbb P^N_T$ 
which satisfy the Chow-field condition \cite[Ch. I, Definition 4.11]{Kol96}. 
Then $Chow_{r, d}(\mathbb P^N/\Z)$ is coarsely represented by 
a semi-normal scheme $\Chow_{r, d}(\mathbb P^N/\Z)$ projective over $\Z$. 
\end{dfn}

\begin{rem}
Since we only need the case when $T$ is a normal noetherian scheme (except for $\Chow_{r, d}(\mathbb P^N/\Z)$), 
let us recall terminologies for this case. 
\begin{enumerate}
\item 
In this case, any well-defined family $U \to T$ of algebraic cycles of $\mathbb P^N/\Z$ 
satisfies the Chow-field condition \cite[Ch. I, Corollary 4.10]{Kol96}. 
\item 
Furthermore, 
if $U = \sum_i m_i U_i$ is a pure $r$-dimensional algebraic cycle 
such that each $U_i$ is flat over $T$, 
then $U \to T$ is a well-defined algebraic families of nonnegative cycles of $\mathbb P^N_T$ 
\cite[Ch. I, Definition 3.10, Definition 3.11, Theorem 3.17]{Kol96}. 
\item 
By construction, 
$\Chow_{r, d}(\mathbb P^N/\Z)$ is the semi normalisation of 
$\Chow'_{r, d}(\mathbb P^N/\Z)$ 
\cite[Ch. I, Definition 3.25.3]{Kol96}, 
where $\Chow'_{r, d}(\mathbb P^N/\Z)$ is a reduced closed subscheme  
of the fine moduli space that parameterises suitable effective divisors, i.e. the projective space corresponding to a linear system. 
Then, 
by \cite[Ch. I, Corollary 3.24.5]{Kol96}, 
the locus $\Chow_{r, d}^{{\rm int}}(\mathbb P^N/\Z)$ parameterising geometrically integral 
cycles is an open subset of $\Chow_{r, d}(\mathbb P^N/\Z)$. 
\end{enumerate}
\end{rem}

\begin{prop}\label{p-chow-var}
Fix positive integers $d$ and $r$. 
Then there exists a flat projective morphism $\pi:V \to S$ 
of quasi-projective $\Z$-schemes 
which satisfies the following property: 
if 
\begin{enumerate}
\item $k$ is a field, 
\item $X$ is an $r$-dimensional geometrically integral projective scheme over $k$, 
and 
\item there is a closed immersion $j:X \hookrightarrow \mathbb P^M_k$ 
over $k$ for some $M \in \Z_{>0}$ such that $(j^*\MO(1))^r \leq d$, 
\end{enumerate}
then there exists a cartesian diagram of schemes: 
\[
\begin{CD}
X @>>> V\\
@VVV @VV\pi V\\
\Spec\,k @>>> S, 
\end{CD}
\]
where the vertical arrows are the induced morphisms. 
\end{prop}

\begin{proof}
We first prove that we may replace the conditions (1)--(3) by the following conditions (1)'--(3)': 
\begin{enumerate}
\item[(1)'] $k$ is an algebraically closed field, 
\item[(2)'] $X$ is an $r$-dimensional projective variety over $k$, 
and 
\item[(3)'] there is a closed immersion $j:X \hookrightarrow \mathbb P^{d+r-1}_k$ 
over $k$ such that $(j^*\MO(1))^r = d$.
\end{enumerate}
Take a triple $(k, X, j:X \hookrightarrow \mathbb P^M_k)$ 
satisfying (1)--(3). 
Note that the claim is equivalent to saying that there are finitely many possibilities for the Hilbert polynomial $\chi(X, j^*\MO(t)) \in \Z[t]$. 
Therefore, passing to the algebraic closure of $k$, 
we may assume that (1)' holds. 
Then (2) and (2)' are equivalent. 
Finally, it follows from \cite[Proposition 0]{EH87} or \cite[Ch. I, Exercise 7.7]{Har77} that 
either $X$ is a projective space or a closed immersion $j:X \hookrightarrow \mathbb P^{d+r-1}_k$. 
We may exclude the former case, thus the problem is reduced to the case when (3)' holds.

Set $N:=d+r-1$ and 
\[
H_1:=\coprod_{\varphi \in \Phi_{r, d}} 
\Hilb_{\mathbb P^N/\Z}^{\varphi} \subset \Hilb_{\mathbb P^N/\Z}, 
\]
where $\Phi_{r, d}$ is the set of polynomials such that 
$\varphi \in \Phi_{r, d}$ if and only if 
there exists an algebraically closed field $k$ and 
a closed immersion $j: X\hookrightarrow \mathbb P^N_k$ 
from an $r$-dimensional projective variety $X$ over $k$ 
such that $(j^*\MO(1))^r = d$. 
Although we do not know yet whether $\Phi_{r, d}$ is a finite set, 
each $\Hilb_{\mathbb P^N/\Z}^{\varphi}$ is a projective $\Z$-scheme. 
For the universal closed subscheme $\Univ_{\mathbb P^N/\Z} \subset \Hilb_{\mathbb P^N/\Z} \times_{\Z} \mathbb P^N_{\mathbb Z}$, 
set $U_1:= \Univ_{\mathbb P^N/\Z} \times_{\Hilb_{\mathbb P^N/\Z}} H_1$. 
In particular, the induced morphism $\rho_1:U_1 \to H_1$ is flat and projective. 
We then define $H_2$ as the open subset of $H_1$ 
such that, for any point $q \in H_1$, 
it holds that $q \in H_2$ if and only if the scheme-theoretic fibre $\rho_1^{-1}(q)$ is 
geometrically integral. 
Let $\rho_2: U_2= U_1 \times_{H_1} H_2 \to H_2$ be the induced flat projective morphism. 
Let $H_3 \to H_2$ be the normalisation of the reduced structure $(H_2)_{\red}$, 
which is a finite morphism. 
Since $H_3$ is normal and $\Chow_{r, d}(\mathbb P^N/\Z)$ is a coarse moduli space, 
the family $U_3=U_2 \times_{H_2} H_3 \to H_3$ induces 
a morphism $\theta:H_3 \to \Chow_{r, d}^{{\rm int}}(\mathbb P^N/\Z)$. 
For any algebraically closed field $k$, the induced map 
$\theta(k):H_3(k) \to \Chow_{r, d}^{{\rm int}}(\mathbb P^N/\Z)(k)$ 
is surjective and any fibre of $\theta(k)$ is a finite set. 
Then, by noetherian induction, 
$H_3$ is of finite type over $\Z$, i.e. $\Phi_{r, d}$ is a finite set. 
Set $\pi:V \to S$ to be $U_3 \to H_3$. 
Then the claim holds.  
\end{proof}


\subsection{Boundedness of regular del Pezzo surfaces}\label{ss2-bounded}

In this subsection, we establish the boundedness 
of geometrically integral regular del Pezzo surfaces (Theorem \ref{t-dP-bdd3}). 
As a consequence, we give a non-explicit upper bound for 
the irregularity  $h^1(X, \MO_X)$ 
(Corollary \ref{c-dP-bdd3}).

\begin{thm}\label{t-dP-bdd1}
Fix a non-negative integer $\epsilon$. 
Then there exists a positive integer 
$d:=d(\epsilon)$ 
which satisfies the following property: 
if $k$ is a field of characteristic $p>0$ and $X$ is a regular del Pezzo surface 
such that $H^0(X, \MO_X)=k$ and $\epsilon(X/k) \leq \epsilon$, 
then there exist a positive integer $N$ and a closed immersion $j:X \hookrightarrow \mathbb P^N_k$ 
such that the degree $(j^*\MO_{\mathbb P^N}(1))^2$ of $j(X)$ is at most $d$. 
\end{thm}

\begin{proof}
By Theorem \ref{t-dP-va2}, $|-12K_X|$ is very ample over $k$. 
Then the assertion follows from Corollary \ref{c-bd-not-gn2}. 
\end{proof}





\begin{thm}\label{t-dP-bdd3}
There exists a flat projective morphism $\rho:V \to S$ 
of quasi-projective $\Z$-schemes 
which satisfies the following property: 
if $k$ is a field of characteristic $p>0$ and 
$X$ is a regular del Pezzo surface over $k$ 
such that $H^0(X, \MO_X)=k$ and $X$ is geometrically reduced over $k$, 
then there exists a cartesian diagram of schemes: 
\[
\begin{CD}
X @>>> V\\
@VV\alpha V @VV\rho V\\
\Spec\,k @>>> S, 
\end{CD}
\]
where $\alpha$ denotes the induced morphism. 
\end{thm}

\begin{proof}
The assertion follows from Proposition \ref{p-chow-var} 
and Theorem \ref{t-dP-bdd1}.  
\end{proof}

\begin{cor}\label{c-dP-bdd3}
There exists a positive integer $h$ 
which satisfies the following property: 
if $k$ is a field of characteristic $p>0$ and 
$X$ is a regular del Pezzo surface over $k$ 
such that $H^0(X, \MO_X)=k$ and $X$ is geometrically reduced over $k$, 
then $\dim_k H^1(X, \MO_X) \leq h$. 
\end{cor}

\begin{proof}
The assertion follows from Theorem \ref{t-dP-bdd3}. 
\end{proof}

\section{Examples}\label{s-example}

In Theorem \ref{t-bd-not-gn}, 
we gave a list of the possibilities 
for the volumes $K_X^2$ of regular del Pezzo surfaces $X$, 
although it depends on $\epsilon(X/k)$. 
Then it is natural to ask whether there actually exists 
a geometrically non-normal example 
which realises each possibility. 
The purpose of this section is to 
give a partial answer by exhibiting several examples. 
We give their construction in Subsection \ref{ss1-example}. 
We then give a summary in Subsection \ref{ss2-example}.

\subsection{Construction}\label{ss1-example}

The purpose of this subsection is 
to construct several regular del Pezzo surfaces which are not geometrically normal.

\begin{ex}\label{e-Fermat-hyp}
Let $\F$ be an algebraically closed field of characteristic $p>0$ 
and let $k:=\F(s_0, s_1, s_2, s_3)$ be 
the purely transcendental extension over $\F$ of degree four. 
Set 
\[
X := \Proj\,k[x_0, x_1, x_2, x_3]/(s_0x_0^p+s_1x_1^p+s_2x_2^p+s_3x_3^p). 
\] 
Then $X$ is a regular projective surface over $k$ such that 
$H^0(X, \MO_X)=k$, 
$(X \times_k \overline k)_{\red} \simeq \mathbb P^2$, and $\epsilon(X/k)=1$ 
\cite[Lemma 9.4, Theorem 9.7]{Tan2}. 
Note that $-K_X$ is ample if and only if $p \in \{2, 3\}$. 
Furthermore, the following hold. 
\begin{enumerate}
\item If $p=2$, then $K_X^2=8$.  
\item If $p=3$, then $K_X^2=3$.  
\end{enumerate}
\end{ex}

\begin{ex}\label{e-Fermat-CI}
Let $\F$ be an algebraically closed field of characteristic two. 
Let 
\[
k:=\F(\{s_i\,|\, 0\leq i \leq 4\} \cup \{t_i\,|\,0 \leq i \leq 4\})
\] 
be the purely transcendental extension over $\F$ of degree ten. 
Set 
\[
X := \Proj\,\left(\frac{k[x_0, x_1, x_2, x_3, x_4]}
{\left(\sum_{i=0}^4 s_ix_i^2, \sum_{i=0}^4 t_ix_i^2\right)}\right). 
\] 
Then $X$ is a regular projective surface over $k$ such that 
$H^0(X, \MO_X)=k$, 
$(X \times_k \overline k)_{\red} \simeq \mathbb P^2$, and $\epsilon(X/k)=2$ 
\cite[Lemma 9.4, Theorem 9.7]{Tan2}. 
Furthermore, we have $K_X^2=4$. 
\end{ex}

\begin{ex}\label{e-geom-quad}
Let $\F$ be an algebraically closed field of characteristic two. 
Let $k:=\F(s_0, s_1, s_2, s_3, s_4)$ be the purely transcendental extension over $\F$ of degree five. 
Set 
\[
X := \Proj\,\left(\frac{k[x_0, x_1, x_2, x_3, x_4]}
{\left(\sum_{i=0}^4 s_ix_i^2, x_0x_1+x_2x_3\right)}\right). 
\] 
Then it holds that $(X \times_k \overline k)_{\red} \simeq \mathbb P^1 \times \mathbb P^1$ and $K_X^2=4$. 
Since $X$ is not geometrically reduced, we have that $\epsilon(X/k)\geq 1$ \cite[Proposition 1.6]{Tan2}. 
Therefore, it follows from Theorem \ref{t-bd-not-gn}(3) that 
$\epsilon(X/k)=1$. 
\end{ex}

\begin{ex}\label{e-Fermat-bup}
Let $\F$ be an algebraically closed field of characteristic two. 
Let $k:=\F(s_0, s_1, s_2, s_3)$ 
be the purely transcendental extension over $\F$ of degree four. 
Set 
\[
Y := \Proj\,k[x_0, x_1, x_2, x_3]/(s_0x_0^2+s_1x_1^2+s_2x_2^2+s_3x_3^2).
\] 
Then $Y$ is a regular projective surface over $k$ such that 
$H^0(Y, \MO_Y)=k$, 
$(Y \times_k \overline k)_{\red} \simeq \mathbb P^2$, 
$\epsilon(Y/k)=1$, and 
$K_Y^2=8$ (Example \ref{e-Fermat-hyp}).
For any $i \in \{0, 1, 2, 3\}$, 
let $C_i$ be the curve on $Y$ defined by $x_i=0$. 
Then we have 
\[
C_3 \simeq \Proj\,k[x_0, x_1, x_2]/(s_0x_0^2+s_1x_1^2+s_2x_2^2). 
\]
The scheme-theoretic intersection $Q:= C_2 \cap C_3$ satisfies 
\[
Q= C_2 \cap C_3 = \Proj\,k[x_0, x_1]/(s_0x_0^2+s_1x_1^2) 
 \simeq \Spec\,k[y]/(y^2+t)
\]
for $t:=s_0/s_1$. 
In particular, $Q$ is a reduced point and $C_2+C_3$ is simple normal crossing. 
Let 
\[
f:X \to Y
\]
be the blowup at $Q$. 
For the proper transform $C'_2$ of $C_2$, 
we have that $(C'_2)^2 = C_2^2 -2 = 0$. 
Since $-K_{C'_2}$ is ample, it holds that $K_X \cdot C'_2<0$. 
Hence, Kleimann's criterion for ampleness 
implies that $X$ is a regular del Pezzo surface. 
Since $(X \times_k \overline k)_{\red}^N$ has a birational morphism to $(Y \times_k \overline k)_{\red}^N \simeq \mathbb P^2$, 
it follows from Theorem \ref{t-bd-not-gn}(3)
$(X \times_k \overline k)_{\red}^N \simeq \mathbb P_{\mathbb P^1}(\MO \oplus \MO(1))$. 
It holds that $K_X^2=K_Y^2-2=6$ and $\epsilon(X/k)=\epsilon(Y/k)=1$, 
where the latter equation follows from Definition \ref{d-epsilon}. 
\end{ex}

\begin{ex}\label{e-cone-deg2}
Let $\F$ be an algebraically closed field of characteristic two. 
Let $k:=\F(s)$ 
be the purely transcendental extension over $\F$ of degree one. 
Set 
\[
Y := \Proj\,k[x, y, z, w]/(x^2+sy^2+zw). 
\] 
It holds that $Y$ is a regular projective surface such that 
$H^0(Y, \MO_Y)=k$, $K_Y^2=8$, and 
$Y \times_k \overline k \simeq \mathbb P(1, 1, 2)$. 
Set 
\[
Y':= D_+(y) \simeq \Spec\,k[x, z, w]/(x^2+s+zw). 
\]
Let $Q \in Y'$ be the closed point defined by 
the maximal ideal $\m:=(x^2+s, z, w)$ of $k[x, z, w]/(x^2+s+zw)$. 
Let $f:X \to Y$ be the blowup at $Q$. 
For the curve $C$ on $Y$ defined by $z=0$, 
we have $C^2=2$. 
Then its proper transform $C'$ satisfies $C'^2= C^2-2= 0$. 
By Kleimann's criterion, we have that $-K_X$ is ample. 
To summarise, we have $K_X^2=K_Y^2-2=6$ and 
$(X \times_k \overline k)^N \simeq \mathbb P(\MO \oplus \MO(2))$, 
where the latter one follows from 
Theorem \ref{t-bd-not-gn}(3) and 
$Y \times_k \overline k \simeq \mathbb P(1, 1, 2)$. 
\end{ex}

\begin{ex}\label{e-Maddock}
Maddock constructed the following examples. 
\begin{enumerate}
\item 
There exists a regular del Pezzo surface $X_1$ 
over a field $k_1$ of characteristic two 
such that $H^0(X_1, \MO_{X_1})=k_1$, 
$X_1$ is geometrically integral over $k_1$, 
$X_1$ is not geometrically normal over $k_1$, 
and $K_{X_1}^2=1$ \cite[Main Theorem]{Mad16}. 
It follows from Theorem \ref{t-bd-not-gn}(3) 
that $(X_1 \times_{k_1} \overline k_1)^N \simeq \mathbb P^2$. 
\item 
There exists a regular del Pezzo surface $X_2$ 
over a field $k_2$ of characteristic two 
such that 
$H^0(X_2, \MO_{X_2})=k_2$, $K_{X_2}^2=2$, 
and $X_2$ is not geometrically reduced. 
It follows from \cite[Proposition 1.6]{Tan2} that $\epsilon(X_2/k_2)\geq 1$. 
By Theorem \ref{t-bd-not-gn}(3), 
it holds that $\epsilon(X_2/k_2)=1$ and 
$(X_2 \times_{k_2} \overline k_2)_{\red}^N \simeq \mathbb P^2$. 
\end{enumerate}
\end{ex}

\subsection{Summary}\label{ss2-example}

We now give a summary of the examples 
established in the previous subsection. 

\begin{table}[H]
\caption{$p=3$ case}
\begin{tabular}{|c|c|c|c|} \hline
$X$ & \hspace{3mm}$K_X^2$\hspace{3mm} & $\epsilon(X/k)$ & 
$(X \times_k \overline{k})_{\red}^N$ \\ \hline \hline
Example \ref{e-Fermat-hyp}  & $3$ & $1$ & $\mathbb P^2$ \\ \cline{2-4} 
\hline 
  \end{tabular}
\end{table}

\begin{table}[H]
\caption{$p=2$ case}
\begin{tabular}{|c|c|c|c|} \hline
$X$ & \hspace{3mm}$K_X^2$\hspace{3mm} & $\epsilon(X/k)$ & 
$(X \times_k \overline{k})_{\red}^N$ \\ \hline \hline
Example \ref{e-Maddock}  & $1$ & $0$ & $\mathbb P^2$ \\ \cline{2-4} 
\hline 
Example \ref{e-Maddock} & $2$ & $1$ &$\mathbb P^2$ \\ \cline{2-4}
\hline
Example \ref{e-Fermat-CI} & $4$ & $2$ &$\mathbb P^2$ \\ \cline{2-4}
\hline
Example \ref{e-geom-quad} & $4$ & $1$ &$\mathbb P^1 \times \mathbb P^1$ \\ \cline{2-4}
\hline
Example \ref{e-Fermat-bup} & $6$ & $1$ &$\mathbb P_{\mathbb P^1}(\MO \oplus \MO(1))$ \\ \cline{2-4}
\hline
Example \ref{e-cone-deg2} & $6$ & $0$ &$\mathbb P_{\mathbb P^1}(\MO \oplus \MO(2))$ \\ \cline{2-4}
\hline
Example \ref{e-Fermat-hyp} & $8$ & $1$ &$\mathbb P^2$ \\ \cline{2-4}
\hline
  \end{tabular}
\end{table}

\begin{rem}
The author does not know whether 
there exists $d \in \Z_{>0}$ such that 
the inequality $K_X^2\leq d$ holds 
for an arbitrary regular del Pezzo surface $X$ 
over a field $k$ with $H^0(X, \MO_X)=k$. 
\end{rem}


\end{document}